\documentclass[12pt]{article}
\usepackage{blindtext}
\usepackage[utf8]{inputenc}

\usepackage{graphics}
\usepackage[
backend=biber,
style=ieee,
sorting=nyt
]{biblatex}
\usepackage{graphicx}
\usepackage{mathtools}
\usepackage{amssymb}
\usepackage{amsfonts}
\usepackage{amsthm}
\usepackage{amsmath}
\usepackage{txfonts}
\usepackage{float}
\usepackage{caption}
\usepackage{subcaption}
\usepackage{comment}
\usepackage{enumitem}
\usepackage{setspace}
\usepackage{hyperref}
\hypersetup{
    colorlinks=true,
    linkcolor=blue,
    filecolor=magenta,      
    urlcolor=cyan,
    }

\newtheorem{theorem}{Theorem}[section]

\newtheorem{lemma}[subsubsection]{Lemma}

\theoremstyle{definition}

\newcommand{\cl}[1]{\mathcal{#1}}
\newcommand{\bb}[1]{\mathbb{#1}}

\newcommand{\G}{\Gamma}

\newcommand{\et}{|\eta|^{1/2}}
\newcommand{\het}{|\hat{\eta}|^{1/2}}
\newcommand{\W}{W(t,\sigma,\eta,\zeta;x_0)}
\newcommand{\heta}{\hat{\eta}}
\newcommand{\hzeta}{\hat{\zeta}}

\DeclareMathOperator{\sgn}{sgn}
\DeclareMathOperator{\supp}{supp}

\allowdisplaybreaks

\title{On some stochastic hyperbolic  equations with symplectic characteristics}
\author{ Enrico Bernardi
  \thanks{
    Dipartimento di Scienze Statistiche
    Paolo Fortunati, Università di Bologna, Bologna,
    Italy. \textbf{e-mail}: enrico.bernardi@unibo.it
  }
  \and
  Leonardo Marconi
  \thanks{
    Dipartimento di Scienze Statistiche Paolo
  Fortunati, Università di Bologna, Bologna, Italy. \textbf{e-mail}:
  leonardo.marconi5@unibo.it}
}

\date{\today}

\begin{document}

\maketitle

\bigskip

\begin{abstract}
		We study the effect of Gaussian perturbations on
                a class of model  hyperbolic partial differential
                equations with double symplectic characteristics in
                low spatial dimensions, extending some recent work
                in \cite{bernardi2023class}. The coefficients of our
                partial differential operators contain harmonic oscillators
                in the space variables, while the noise is additive,
                white in time and colored in space.
                We provide  sufficient conditions on the spectral
                measure of the covariance functional describing the
                noise that allows for the existence of a random field
                solution for the resulting stochastic partial
                differential equation.
                Furthermore we show how the symplectic structure of the set of
                multiple points affects the regularity of the noise
                needed to build a measurable process solution.
                Our approach is based on some explicit computations for the
                fundamental solutions of several model partial
                differential operators together with their explicit Fourier transforms.    
\end{abstract}

	
	
\allowdisplaybreaks



\section{Introduction}
A number of recent papers, tapping into the powerful techniques
presented in the seminal works
\cite{Walsh1986AnIT}, \cite{dalang1999extending} and
\cite{dalang1998stochastic}, have extended several  classical deterministic
results for solutions of hyperbolic operators of various types -linear, semilinear,
and whose principal symbols have multiple involutive or symplectic characteristics- to the stochastic
framework:
for a current review one can consult, among others,  \cite{ascanelli2019solution},
\cite{abdeljawad2020deterministic}, \cite{ascanelli2020random},
\cite{ascanelli2018random}, \cite{bernardi2023class}.
Our goal in this work is to continue the study of possible
extensions to the stochastic framework by analyzing three model linear
stochastic operators, we deem not having been previously examined in the literature.
In particular, we would like to understand how the symplectic geometry of
the principal symbols and the conditions on the lower order terms
influence the construction of the random field solutions via the
corresponding colors of the noise, an analysis started in \cite{bernardi2023class}.

 The general problem is presented here in a compact form:
\begin{equation}
\label{eq}
    \begin{cases}
        P_iu_i=\Dot{F}_i\left(t,x\right)& t\geq 0, x \in \bb{R}^i,i \in \{1,2,3\}\\
        u_i\left(0,x,y,z\right)=0\\\partial_t u_i\left(0,x,y,z\right)=0,
    \end{cases}
\end{equation} 
where formally
\[F_i(\phi):= \int_{\bb{R}^{i+1}}\phi(t,x)\Dot{F}_i(t,x)\,dt dx, \quad \phi \in C_0^\infty\left(\bb{R}^{i+1}\right),\]
is a mean zero family of normal random variables defined on a probability space $(\Omega,\cl{F},\bb{P})$ with covariance
\begin{equation}
	\label{covariance}
	E(F_i(\phi)F_i(\psi))=\int_0^\infty\,dt\int_{\bb{R}}\,dx \int_{\bb{R}}\,dy\, \phi(t,x)f_i(x-y)\psi(t,y),
\end{equation}
for $\phi,\psi \in C_0^\infty\left(\bb{R}^{i+1}\right)$.     
The operators $P_i,i \in \{1,2,3\}$ are defined as
\begin{gather}
	\label{one dim operator}
	P_1:=D_t^2-(D_x^2+x^2),\\
	\label{two dim operator}
	P_2:=D_t^2-\mu(D_x^2+x^2D_y^2)+bD_y, \quad |b|<\mu\\
	\label{three dim operator}
	P_3:=D_t^2-\mu(D_x^2+x^2D_y^2)-aD_z^2+bD_y, \quad |b|<\mu,\, a >0.
\end{gather}
$P_1$ is a $0$-th order perturbation of the wave operator, which is
strictly hyperbolic, and we use it mostly as an introduction to some
of the techniques used in the latter cases.
Regarding the two- and three-dimensional operators $P_2$ and $P_3$, their corresponding principal symbols are
\begin{gather*}
	\label{two dim symbol}
	p_2:=\tau^2-\mu(\xi^2+x^2\eta^2),\\
	\label{three dim symbol}
	p_3:=\tau^2-\mu(\xi^2+x^2\eta^2)-a\zeta^2,
\end{gather*}
which are readily seen to be  hyperbolic w.r.t. $\tau$ and vanishing at the second order on the $C^\infty$-manifolds
\begin{gather*}
	\label{two dim manifold}
	\Sigma_2:=\{(t,x,y,\tau,\xi,\eta)\in \dot{T}^*\bb{R}^3,x=\tau=\xi=0\},\\
	\label{three dim manifold}
	\Sigma_3:=\{(t,x,y,z,\tau,\xi,\eta,\zeta)\in \dot{T}^*\bb{R}^4,x=\tau=\xi=\zeta=0\},
\end{gather*}
where $\dot{T}^*\bb{R}^i$ denotes the phase space cotangent bundle minus the
zero section. The manifold $\Sigma_2$ is purely symplectic in its
space variables and covariables, i.e. the two-form $ d\xi \wedge dx  $
is non degenerate or equivalently the Poisson brackets of the
functions defining the manifold do not all vanish on the manifold
itself, while $\Sigma_3$ presents a mixed
involutive-symplectic structure. Some review of the classical
symplectic framework for differential equations with multiple
characteristics can be found, among others, in \cite{bernardi2011cauchy} ,
\cite{bernardi2023class} and \cite{bony1991quadratic}.

Calling $H_i$ the Hamiltonian field of $p_i$, the fundamental matrix $F_i(\omega)$ associated  to $p_i$ at a point $\omega \in \Sigma_i$ is computed as $F_i(\omega)=\frac{1}{2}dH_i(\omega)$.
One can easily check that for all points $\omega_i \in \Sigma_i$ both
$F_2(\omega_2)$ and $F_3(\omega_3)$ have just two non-zero complex
eigenvalues $\pm i \mu \eta$.
One then defines in this case the corresponding operator to be of non-effectively hyperbolic type, \cite{bony1991quadratic}.
The request of well-posedness of the Cauchy problem in the $
C^{\infty}$ and Gevrey classes for this type of operators gives rise to the (strict) Levi type
condition $|b|<\mu$, where $\mu$ is what is usually known as the
positive trace of the associated harmonic oscillator (see e.g. \cite{hormander1977cauchy} and \cite{bony1991quadratic}).
The (non-strict) Levi condition $ |b| \leq \mu $ is, as it is well known,
necessary -and in its strict version
sufficient as well- in order to reach
the well-posedness of the Cauchy problem for $ P_{3} $ in the $
C^{\infty} $ category, see e.g. \cite{ivrii1974necessary}. It will be
shown below how, even in the stochastic construction, this condition on
the lower order terms plays a fundamental role for the explicit existence
of the random field solution.

The hyperbolic operator with symbol $ p_{2} $ has been studied in the
deterministic context in a number of papers, see
for instance \cite{melrose1984wave}, \cite{lascarB},\cite{lascarBlascarR}, aiming to calculate Poisson formulas
for hypoelliptic operators on compact manifolds and studying the
propagation of singularities for their solutions.
Both these objectives rely on having at one's disposal a 
rather explicit expression for the fundamental solution; thus in order to develop
similar results in a stochastic setup we must necessarily possess a rather sizeable
knowledge of such a fundamental solution. That  being said, the presence of harmonic oscillators in the symbols $ p_{2} $ and $ p_{3} $ together with the symplectic nature of
their double manifolds makes the
more traditional approaches of extending the classical calculi of
pseudo-differential analysis to the stochastic setting much more
complicated. This is essentially due to the fact that the simple
bicharacteristic curves exhibit a periodic nature, making it very hard
to present any fundamental solution as an integral
Fourier operator whose phase is supported over those very curves. An
explicit procedure is then required and this is what will be done in the present work
through suitable Hermite functions expansions and a precise control
over the related coefficients.

Finally, we observe that the arguments we utilize in our analysis
could very well be used to deal with a more general second order
hyperbolic quadratic form with symplectic characteristics.
Nevertheless, the simpler nature of our model
examples allows us to remove some of the complexities in
\cite{melrose1984wave} or \cite{lascarB}, where only the
codimension $ 3 $ was considered, and develop a sharper understanding of the structure of
our fundamental solution processes. We observe that the symplectic nature of the
multiple manifold is not by itself an  a-priori  obstacle to obtaining a
closed form fundamental solution, as was proven in
\cite{bernardi2023class}. What makes the problem more challenging in our case
is indeed the presence of the harmonic oscillators.


The function $f_i:\bb{R}^{i}\to \bb{R}$ is taken to be continuous except at most at zero, and even. This requirements are necessary to guarantee that the functional expressed in (\ref{covariance}) is non-negative definite. As seen in \cite{schwartz2002theorie} ,  this is also equivalent to the existence of a non-negative tempered measure $\nu_i$ whose Fourier transform is $f_i$, i.e., for all $\phi \in\cl{S}(\bb{R}^i)$
\begin{equation*}
\int_{\bb{R}^i}\phi(x)f_i(x)dx=\int_{\bb{R}^i}\cl{F}\phi(\xi)\nu_i(d\xi),
\end{equation*}
where, and henceforth, $\cl{F}\phi$ denotes the Fourier transform
\begin{equation*}
    \cl{F}\phi(\xi):=\int_{\bb{R}^d}e^{-ix\cdot \xi}\phi(x)dx
\end{equation*}
and
\begin{equation*}
    \cl{F}^{-1}\phi(x):=\frac{1}{2\pi}\int_{\bb{R}^d}e^{ix\cdot \xi}\cl{F}\phi(\xi)d\xi.
\end{equation*}
Denoting with $E_i$ the fundamental solution of the operator $P_i$, we will say that $u_i(t,x)$ defined as
\begin{equation}
	\label{formal solution}
	u_i(t,x)=\int_0^t \int_{\bb{R}^d}E_i\left(t-s,x-y\right)\Dot{F}_i(s,y)ds\,dy
\end{equation} 
is a \textit{random field solution} to (\ref{eq}) if  the stochastic integral is well defined and the map $(t,x)\to u_i(t,x)$ is measurable. For the theory of stochastic integration, we refer to \cite{dalang1999extending}, where the author performs an extension and adaptation of Walsh's construction of the martingale measure stochastic integral (see \cite{Walsh1986AnIT}). We will shortly describe it.\\
We denote with $\cl{D}\left(\bb{R}^p\right)$ the space of functions $\phi \in C_0^\infty\left(\bb{R}^p\right)$ endowed with the topology described by the following notion of convergence. Given a sequence $(\phi_n)_{n \in \bb{N}} \subset \cl{D}\left(\bb{R}^p\right)$ and a function $\phi \in \cl{D}\left(\bb{R}^p\right)$ we say that $\phi_n$ converges to $\phi$ and write $\phi_n \to \phi$ if:
\begin{itemize}
	\item there exists a compact set $K \subset \bb{R}^p$ such that $\supp\left(\phi_n-\phi\right)\subset K$ for all $n \geq 1$,
	\item $\lim_{n \to \infty}D^\alpha \phi_n =D^\alpha \phi$ uniformly in $K$ for every multi-index $\alpha$.
\end{itemize}
Let $\cl{B}_b(\bb{R}^i)$ e the $\sigma$-field of all bounded Borel sets of $\bb{R}^i$. The first step is to extend $F_i$ to a \textit{worthy martingale measure} (see \cite{Walsh1986AnIT}). Using (\ref{covariance}) we can verify that $F_i$ is $L^2$-continuous and, as such, we can extend it to a $\sigma$-finite $L^2$-valued measure by approximating indicator functions of sets in $\cl{B}_b(\bb{R}_+\times \bb{R}^i)$ with elements of $\cl{D}\left(\bb{R}^{i+1}\right)$. We set
\begin{equation*}
    M_{i,t}(B)=F_i([0,t]\times B), \quad B \in \cl{B}_b(\bb{R}^i)
\end{equation*}
and 
\begin{equation*}
    \cl{F}_{i,t}^0=\sigma \left(M_{i,s}(B), s \leq t , B \in \cl{B}_b(\bb{R}^i)\right), \quad \cl{F}_{i,t}=\cl{F}_{i,t}^0 \vee \cl{N},
\end{equation*}
with $\cl{N}$ being the $\sigma$-field generated by the $\bb{P}$-null sets. We then have that, by construction, $t \to M_t(B)$ is a continuous martingale and
\begin{equation*}
	F_i(\phi)=\int_{\bb{R}^+}\int_{\bb{R}^i}\phi(t,x)M(dt,dx), \quad \phi \in \cl{D}(\bb{R}^{i+1})
\end{equation*} 
We call a function $(s,x,\omega) \to g(s,x,\omega)$ \textit{elementary} if it is of the form
\begin{equation*}
    g(s,x,\omega)=1_{(a,b]}(s)1_A(x)X(\omega) a,b \in \bb{R},\ 0\leq a<b, \ A \in \cl{B}_b(\bb{R}^d), 
\end{equation*}
where $a,b \in \bb{R}, \ a<b, \ A \in \cl{B}_b(\bb{R}^i) $ and $X$ is a $\cl{F}_{i,a}$-measurable r.v.
We denote via $\cl{E}$ the space of all finite linear combinations of elementary functions and call \textit{predictable} the $\sigma$-field on $\bb{R}_+ \times \bb{R}^i\times \Omega$ generated by the elements of $\cl{E}$.
Moreover, set
\begin{equation}
    \label{norm +}
    ||g||_+:=E\left(\int_0^t ds \int_{\bb{R}^i} dy\int_{\bb{R}^i} dx |g(s,x,\cdot)|f_i(x-y)|g(s,y,\cdot)|\right)
\end{equation}
and
\begin{equation}
    \label{norm 0}
    ||g||_0:=E\left(\int_0^t ds \int_{\bb{R}^i} dy\int_{\bb{R}^i} dx g(s,x,\cdot)f_i(x-y)g(s,y,\cdot)\right).
\end{equation}
In \cite{Walsh1986AnIT}, Walsh defines the martingale-measure stochastic integral on the complete space $\cl{P}_+$ of predictable functions $g$  with  $||g||_+< +\infty$.\\
On the other hand, Dalang in his work \cite{dalang1999extending} is able to extend such a construction defining the martingale
\begin{equation*}
    t \to \int_0^t \int_{\bb{R}^i} v(s,x,\cdot) M\left(dx, ds\right)
\end{equation*}
for all elements $v$ of the completion $\cl{P}_0$ of $\left(\cl{E},||\,\cdot\,||_0\right)$. Moreover, denoting as $\Bar{\cl{P}}$ the space of all predictable functions $h(t,x,\omega)$ such that $x \to h(t,x,\omega) \in \cl{S}'(\bb{R}^i)$ for every $(t,\omega) \in [0,T]\times \Omega $, $ \cl{F}h(t,\cdot,\omega)(\xi)$ is a function a.s. and
\begin{equation}
\label{Pbar bound}
    \left|\left|h\right|\right|_0':=E\left(\int_0^t ds \int_{\bb{R}^i} \left|\cl{F}h(t,\cdot,\omega)(\xi)\right|^2 \nu_i(d\xi)  \right)^{1/2}< +\infty,
\end{equation}
and calling $\cl{E}_0$ the subset of $\cl{P}_+ $ consisting of all functions $g(s,x,\omega)$ such that $x \to g(s,x,\omega)\in \cl{S}(\bb{R}^i)$, then we can identify $\cl{P}_0$ with the set of elements $h(t,x,\omega)$ of $\Bar{\cl{P}}$ such that it exists a sequence $h_n(t,x,\omega)$ of elements of $\cl{E}_0$ that gives 
\begin{equation}
\label{sequence in E0}
    \lim_{n\to \infty}||h_n-h||_0'=0.
\end{equation}
Furthermore, we observe that this implies that for elements $h(s,x,\omega) \in\cl{P}_0$
\begin{equation*}
    ||h||_0=||h||_0'.
\end{equation*}
Finally, we define the physicist's Hermite polynomials and Hermite functions as
\begin{gather}
	\label{hermite}
	H_n(x)=(-1)^n e^{x^2} \frac{d^n}{d x^n}e^{-x^2}, \nonumber \\
	\Psi_n(x)=\frac{1}{(2^n n! \sqrt{\pi})^{1/2}}e^{-x^2\over{2}}H_n(x),
\end{gather}
and name
\begin{gather*}
	\rho_n\left(\eta,\zeta\right):=\mu |\eta|\left(2n+1\right)+a\zeta^2 -b\eta\\
	\hat{\rho}_n(\eta):=\rho_n(\eta,0).
      \end{gather*}
We are now ready to state our main result, which contains all the several claims proven below, itemized
according to the increasing number of space variables.
Recalling that our stochastic PDEs are
\begin{equation*}
	\begin{cases}
		P_iu_i=\Dot{F}_i\left(t,x\right)& t\geq 0, x \in \bb{R}^i,i \in \{1,2,3\}\\
		u_i\left(0,x,y,z\right)=0\\\partial_t u_i\left(0,x,y,z\right)=0,
	\end{cases}
\end{equation*} 
where the operators $P_i,i \in \{1,2,3\}$ are defined as
\begin{gather*}
	P_1:=D_t^2-(D_x^2+x^2),\\
	P_2:=D_t^2-\mu(D_x^2+x^2D_y^2)+bD_y, \quad |b|<\mu\\
	P_3:=D_t^2-\mu(D_x^2+x^2D_y^2)-aD_z^2+bD_y, \quad |b|<\mu,\, a >0,
\end{gather*}
we have the following:
\begin{theorem}
	\label{theorem 1}
\begin{enumerate}[label=(\roman*)]
	\item Let $\nu_1(d\xi_0)=d\xi_0$ be the Lebesgue measure. Then the one-dimensional problem has a random field solution given by
	\begin{equation*}
		u_1(t,x)=\int_0^t \int_{\bb{R}^d}E_1\left(t-s,x;x_0\right)\Dot{F}_1(s,x_0)ds\,dx_0.
	\end{equation*}
	with
	\[E_1(t,x,x_0)=-H(t)\sum_{n=0}^\infty \frac{\sin(\sqrt{2n+1}\,t)}{\sqrt{2n+1}}\Psi_n(x_0)\Psi_n(x).\].
	\item Let $\nu_2=(d\xi,d\heta)=d\xi\hat{\nu}_2(d\heta)$ and assume that
	\begin{equation}
		\label{nu2 condition}
		\int_{\bb{R}}\frac{1}{1+|\heta|^{1/2}}\hat{\nu}_2(d\heta) <+\infty;
	\end{equation}
	then the two-dimensional problem has a random field solution given by
	\begin{equation*}
		u_2(t,x,y)=\int_0^t \int_{\bb{R}^d}E_2\left(t-s,x,y-y_0;x_0\right)\Dot{F}_1(s,y_0,x_0)ds\,dx_0\,dy_0.
	\end{equation*}
	with
		\begin{align*}
		E_2&\left(t,x,y;x_0\right)=-\frac{H\left(t\right)}{2\pi}\\\times&\int_{\bb{R}}e^{iy\eta}\et\sum_{n=0}^\infty\frac{\sin\left(\hat{\rho}_n\left(\eta\right)^{1/2} t\right)}{\hat{\rho}_n\left(\eta\right)^{1/2}}\Psi_n\left(x_0 \et\right) \Psi_n\left(x \et\right) d\eta.
	\end{align*}
	\item Let $\nu_3(d\xi,d\heta,d\hzeta)=d\xi \hat{\nu}_3(d\heta,d\hzeta)$. Assume $\hat{\nu}_3$ is absolutely continuous w.r.t. the Lebesgue measure and that it admits a density of the form $(\heta,\hzeta)\to w\left(|\heta|^2+|\hzeta|^2\right)$ for which it exists an $\alpha<1/3$ such that
	\begin{equation}
		\label{nu3 condition}
		\int_{\bb{R}^2}\frac{1}{\left(1+|\heta|^2+|\hzeta|^2\right)^{\alpha}}\hat{\nu}_3(d\heta,d\hzeta) <+\infty.
	\end{equation}
	Then the three-dimensional problem has a random field solution given by
	\begin{align*}
		u_3(t,x,y,z)=\int_0^t \int_{\bb{R}^d}&E_3\left(t-s,x,y-y_0,z-z_0;x_0\right)\\ \times&\Dot{F}_3(s,x_0,y_0,z_0)ds\,dx_0\,dy_0\,dz_0.
	\end{align*}
	with
	\begin{align*}
		E_3&\left(t,x,y,z;x_0\right)=-\frac{H\left(t\right)}{2\pi}\\\times&\int_{\bb{R}^2}e^{iy\eta+iz\zeta}\et\sum_{n=0}^\infty\frac{\sin\left(\rho_n\left(\eta,\zeta\right)^{1/2} t\right)}{\rho_n\left(\eta,\zeta\right)^{1/2}}\Psi_n\left(x_0 \et\right) \Psi_n\left(x \et\right) d\eta\, d\zeta.
	\end{align*}
\end{enumerate}
\end{theorem}

\vskip 24pt
We would like to highlight again how the color of the
noise, observable through the corresponding measures, relates to the geometry of
the double manifolds.

The plan of the paper is as follows:
in section \ref{one dim} we  organize the proof for the
one-dimensional case, parts of which are to be used in the later
sections. In particular in subsection \ref{fd1} the formal fundamental
solutions are computed and in subsection \ref{rfs1} the corresponding
random fields are constructed.
Section \ref{23dims} contains the main arguments. In subsection
\ref{fund sol for two and three dim} the formal fundamental solutions
are explicitly calculated for both codimension $ 3 $ and $ 4 $. Subsections
\ref{two dim} and \ref{3dib} are devoted to the estimate of the formal
solutions obtained, and in the final subsection \ref{fmeas} the
existence of the random field solutions as stated in Theorem
\ref{theorem 1} is eventually proven.

\section*{Acknowledgements}

The authors would like to thank Professor A. Lanconelli for numerous
helpful discussions and feedback.

\section{The case of one spatial dimension}
\label{one dim}
\subsection{Fundamental solution}
\label{fd1}
We now proceed with the case of one spatial dimension: this, besides
proving the related simple case, helps us collecting all the major
tools needed for the subsequent parts.

We consider the problem
\[P_1E_1(t,x;x_0)=\delta(t)\delta(x-x_0).\]
Performing an Hermite series expansion of $x \to E_1(t,x;x_0)$ and using the well-known identities (see \cite{lebedev1965special} or \cite{szeg1939orthogonal})
    \begin{equation}
    \label{mehler}
        \sum_{n=0}^\infty\Psi_n(x)\Psi_n(x_0) = \delta(x-x_0),
    \end{equation}
    
    \begin{equation}
    \label{hermite differential identity}
        (\partial_x^2-x^2+2n+1)\Psi_n(x)=0 ,
    \end{equation}
we get
\begin{equation}
\label{hermite decomposition}
   P_1 \sum_{n=0}^\infty E_{1,n}(t,x_0) \Psi_n(x)=\sum_{n=0}^\infty (P_{1,n} E_{1,n})(t,x_0) \Psi_n(x)=\sum_{n=0}^\infty \delta(t)\Psi_n(x)\Psi_n(x_0),
\end{equation}
where 
\[P_{1,n}=D_t^2-(2n+1).\]
Equating term by term in (\ref{hermite decomposition}) we obtain the ODE problem
\[\begin{cases}
    P_{1,n} E_{1,n}(t,x_0)=\delta(t)\Psi_n(x_0),\\
    E_{1,n}(0)=\partial_tE_{1,n}(0)=0,
\end{cases}\]
whose solution is
\begin{align}
\label{ode}
    E_{1,n}(t,x_0)=&-\Psi_n(x_0)\int_0^t \delta(t-\zeta)\frac{\sin(\sqrt{2n+1}\,\zeta)}{\sqrt{2n+1}}\,d\zeta \nonumber \\&=-H(t)\frac{\sin(\sqrt{2n+1}\,t)}{\sqrt{2n+1}}\Psi_n(x_0),
\end{align}
where $H\left(t\right)$ denotes the Heavyside function,i.e,
\[H\left(t\right)=\begin{cases}
    1, &t \geq 0,\\
    0, &\textrm{otherwise.}
\end{cases}\]
Hence
\[E_1(t,x;x_0)=-H(t)\sum_{n=0}^\infty \frac{\sin(\sqrt{2n+1}\,t)}{\sqrt{2n+1}}\Psi_n(x_0)\Psi_n(x).\]
Now, the Hermite functions are eigenvectors of the Fourier transform (see e.g.\cite{szeg1939orthogonal}); in particular
\begin{equation*}
    \cl{F}\Psi_n(\eta)=(-i)^n\Psi_n(\eta)
\end{equation*}
and so
\begin{align}
\label{fund sol 1 dim}
    \mathtt{F}_1(t,x;\xi_0):=&\cl{F} E_1(t,x;\cdot)\left(\xi_0\right)\nonumber\\=&-H(t)\sum_{n=0}^\infty (-i)^n\frac{\sin\left(t\sqrt{2n+1}\right)}{\sqrt{2n+1}}\Psi_n(\xi_0)\Psi_n(x)
\end{align}
\subsection{Random field solution}
\label{rfs1}
For the reader's convenience, we recall the expression for the candidate random field solution
\begin{equation*}
	u_1(t,x)=\int_0^t \int_{\bb{R}^d}E_1\left(t-s,x;x_0\right)\Dot{F}_1(s,x_0)ds\,dx_0.
\end{equation*}
In order to prove that the expression above is indeed a real-valued process, we need to show that $E_1\left(t-s,x;x_0\right)$ is an element of the space $\cl{P}_0$. To do that, the first step is to show the following.
\begin{lemma}
\label{lemma 1 1dim}
    Let $\mathtt{F}_1(t,x;\xi_0)$ be defined as in (\ref{fund sol 1 dim}). Then
    \begin{equation*}
    I:=\int_0^t \int_{\bb{R}}|\mathtt{F}_1\left(t-s,x;\xi_0\right)|^2 d\xi_0\,ds <+\infty.
\end{equation*}
\end{lemma}

\begin{proof}
By Parseval identity, we have
\begin{align*}
    I=\int_0^t\sum_{n=0}^{+\infty} \frac{\sin^2\left((t-s)\sqrt{2n+1}\right)}{2n+1}\Psi_n^2(x)ds\leq t\sum_{n=0}^{+\infty}\frac{1}{2n+1}\Psi_n^2(x)
\end{align*}
Now, we recall that there exist constants $C$ and $D$ such that, (see e.g. \cite{bonan1990estimates})
\begin{equation}
\label{bound on psi2}
    Cn^{-1/6}<\max_{x\in\bb{R}}\Psi_n^2(x)<Dn^{-1/6}\quad n\geq 1,
\end{equation}
This implies that
\begin{equation*}
    I \leq Dt\left(\Psi_0^2(x)+\sum_{n=1}^{+\infty}(2n+1)^{-1}n^{-1/6}\right)\leq Dt\left(\pi^{-1/4}+\sum_{n=1}^{+\infty}(2n+1)^{-1}n^{-1/6}\right) 
\end{equation*}
which is bounded.
\end{proof}
This shows that the fundamental solution $E_1\left(t-s,x;x_0\right)$ belongs to the space $\Bar{\cl{P}}$. However, to prove its membership to $\cl{P}_0$ as well we need to find a sequence of functions $\left(E_{1_m}\right)_{m \in \bb{N}}$  in $\cl{E}_0$ converging to $E_1$ in $\Bar{\cl{P}}$.\\
Unsurprisingly, we set
\begin{align*}
	E_{1,m}\left(t-s,x;x_0\right):=-H(t)\sum_{n=0}^m \frac{\sin\left((t-s)\sqrt{2n+1}\right)}{\sqrt{2n+1}}\Psi_n(x_0)\Psi_n(x).
\end{align*}
Firstly, we observe that for each $m \in \bb{N}$, $x_0 \to E_{1,m}\left(t,x;x_0\right) \in \cl{S}\left(\bb{R}\right)$.
Moreover, setting
\begin{align*}
	\mathtt{F}_1^m\left(t-s,x;\xi_0\right):=-H(t)\sum_{n=m}^\infty (-i)^n\frac{\sin\left((t-s)\sqrt{2n+1}\right)}{\sqrt{2n+1}}\Psi_n(\xi_0)\Psi_n(x),
\end{align*}
we have
\begin{align*}
	\left|\left|E_{1_m}-E\right|\right|_0=&\int_0^t ds\int_{\bb{R}}\left|\mathtt{F}_1^m\left(t-s,x;\xi_0\right)\right|^2d\xi_0.
\end{align*}
Hence, from the proof of Proposition \ref{lemma 1 1dim}, we gather
\begin{align*}
	\left|\left|E_{1_m}-E\right|\right|_0 \leq Dt\sum_{n=m}^{+\infty}(2n+1)^{-1}n^{-1/6}\xrightarrow{m \to +\infty}0.
\end{align*}
Now, it remains to prove measurability, and we do that by checking $\bb{L}^2(\Omega)$-continuity of the solution process in the separate variables. 
Denoting
\begin{equation*}
	S_1(t,x,h)=\sum_{n=0}^{+\infty} \frac{\sin^2\left(t\sqrt{2n+1}\right)}{2n+1}\left(\Psi_n(x+h)-\Psi_n(x)\right)^2,
\end{equation*}
the increment in space gives
\begin{align*}
	&\bb{E}\left[\left|u_1(t,x+h)-u_1(t,x)\right|^2\right]\\=&\int_0^tds\int_{\bb{R}} \left|\mathtt{F}_1\left(t-s,x+h;\xi_0\right)- \mathtt{F}_1\left(t-s,x;\xi_0\right)\right|^2d\xi_0\\=&\int_0^tS_1(t-s,x,h)ds.
\end{align*}
Thanks to (\ref{bound on psi2}) we have
\begin{equation*}
	S_1(t,x,h)\leq 4D\left(\pi^{-1/4}+\sum_{n=1}^{+\infty}(2n+1)^{-1}n^{-1/6}\right)<+\infty,
\end{equation*}
granting uniform convergence and consequently the continuity of $S_1(t-s,x,h)$ w.r.t. his last variable.
Therefore, by dominated convergence
\begin{equation}
	\lim_{h\to 0}\bb{E}\left[\left|u_1(t,x+h)-u_1(t,x)\right|^2\right]=0.
\end{equation}
For the increment in time we write
\begin{align*}
	&\bb{E}\left[\left|u_1(t+h,x)-u_1(t,x)\right|^2\right]\\ \leq&2\int_0^tds\int_{\bb{R}} \left|\mathtt{F}_1\left(t+h-s,x;\xi_0\right)- \mathtt{F}_1\left(t-s,x;\xi_0\right)\right|^2d\xi_0 \\\mathrel{\phantom{=}}&+2\int_t^{t+h}ds\int_{\bb{R}} \left|\mathtt{F}_1\left(t+h-s,x;\xi_0\right)\right|^2d\xi_0\\=&2I_1(h)+2I_2(h).
\end{align*}
Lemma \ref{lemma 1 1dim} gives directly
\begin{equation*}
	\lim_{h \to 0}I_2(h)=0. 
\end{equation*}
Similarly to the increment in the $x$ variable, we call
\begin{equation*}
	S_1'(t,x,h):=\sum_{n=0}^{+\infty}\frac{\left(\sin\left((t+h)\sqrt{2n+1}\right)-\sin\left(t\sqrt{2n+1}\right)\right)^2}{\sqrt{2n+1}}\Psi_n^2(x)
\end{equation*}
and by Parseval's identity, we obtain
\begin{align*}
	I_1(h)=\int_0^tS_1'(t-s,x,h)ds.
\end{align*}
For $S_1'(t-s,x,h)$, a bound similar to (\ref{sigma bound}) holds. So, by dominated convergence
\begin{equation*}
	\lim_{h\to 0}I_2(h)=0.
      \end{equation*}

      \section{The case of two and three spatial dimensions}
      \label{23dims}
\subsection{Common computations: fundamental solution}
\label{fund sol for two and three dim}
Here we explicitly derive the fundamental solution of the operator $P_3$. Since the structure of the calculation is much the
same as in the two-dimensional case, we omit the latter and recover the corresponding fundamental solution at the end of the section.

Hence, our problem becomes now
\begin{equation*}
	P_3E_3=\delta(t)\delta(x-x_0)\delta(y)\delta(z).
\end{equation*}
Taking the Fourier transform w.r.t. $\left(y,z\right)$ we get
\begin{equation}
	\label{problem 3d}
	\widehat{P_3E_3}=\left(D_t^2 -\mu\left(D_x^2+ x^2\eta^2\right)-a\zeta^2 +b\eta\right)\hat{E}_3
	\left(t,x,\eta,\zeta;x_0\right)=\delta\left(t\right)\delta\left(x-x_0\right).
\end{equation}
Putting $\sigma= |\eta|^{1/2}x$ we define
\[W\left(t,\sigma,\eta,\zeta;x_0\right)=\hat{E}_3\left(t,|\eta|^{-1/2}x,\eta,\zeta;x_0\right);\]
thus
\[\partial_\sigma^2 W\left(t,\sigma,\eta,\zeta;x_0\right) = |\eta|^{-1} \partial_x^2 \hat{E}_3\left(t,|\eta|^{-1/2}x,\eta,\zeta;x_0\right) \]
and, consequently
\[\left(D_x^2+x^2\eta^2\right)\hat{E}_3\left(t,|\eta|^{-1/2}x,\eta,\zeta;x_0\right)=|\eta|\left(D_\sigma^2+\sigma^2\right) W\left(t,\sigma,\eta,\zeta;x_0\right).\]
Therefore, since Dirac's delta is homogeneous of degree $-1$, naming $\Lambda:=D_t^2 -\mu |\eta|\left(D_\sigma^2+\sigma^2\right)-a\zeta^2 +b\eta$,  (\ref{problem 3d}) finally becomes
\begin{align}
	\label{eq2}
	\Lambda W\left(t,\sigma,\eta,\zeta;x_0\right)
	=\et  \delta\left(t\right)\delta\left(\sigma-x_0\et\right).
\end{align}
We expand $\sigma \to \W$ in Hermite functions 
\[\W=\sum_{n=0}^\infty W_n\left(t,\eta,\zeta;x_0\right) \Psi_n\left(\sigma\right)\]
and define
\[\Lambda_n=D_t^2 -\mu |\eta|\left(2n+1\right)-a\zeta^2 +b\eta
		.\]
		It follows that (\ref{eq2}) transforms into
		\begin{align*}
			\Lambda \W =&\Lambda\sum_{n=0}^\infty W_n\left(t,\eta,\zeta;x_0\right) \Psi_n\left(\sigma\right)\\=&\sum_{n=0}^\infty \Lambda_n W_n\left(t,\eta,\zeta;x_0\right) \Psi_n\left(\sigma\right)\\=&\delta\left(t\right)\sum_{n=0}^\infty\Psi_n\left(x_0\et\right)\Psi_n\left(\sigma\right),
		\end{align*}
		and we split the sum term by term, equating the coefficients of the expansion
		\begin{equation*}
			\left(\Lambda_n W_n\right)\left(t,\eta,\zeta;x_0\right)=\delta\left(t\right)\Psi_n\left(x_0\et\right)\et,
		\end{equation*}
		so that we may focus on the ODE problem
		\begin{equation}
			\label{ode eq}
			\begin{cases}
				\left(\partial_t +\rho_n\left(\eta,\zeta\right)\right)W_n\left(t,\eta,\zeta;x_0\right)=-\delta\left(t\right)\Psi_n\left(x_0 \et\right) \et\\
				W_n\left(0,\eta,\zeta;x_0\right)=\partial_t W_n\left(0,\eta,\zeta;x_0\right)=0,
			\end{cases}
		\end{equation}
		which is solved by
		\begin{align*}
			W_n\left(t,\eta,\zeta;x_0\right)=&-\Psi_n\left(x_0\et\right)\et \cdot \int_0^t \delta\left(t-\tau\right) \frac{\sin\left( \rho_n\left(\eta,\zeta\right) \tau\right)}{ \rho_n\left(\eta,\zeta\right)}
			\\=&-H\left(t\right)\Psi_n\left(x_0 \et\right)\et \frac{\sin\left(\rho_n\left(\eta,\zeta\right)^{1/2} t\right)}{\rho_n\left(\eta,\zeta\right)^{1/2}}.
		\end{align*}
		Therefore, we find
		\begin{equation*}
			\W=  -H\left(t\right)\et\sum_{n=0}^\infty\frac{\sin\left(\rho_n\left(\eta,\zeta\right)^{1/2} t\right)}{\rho_n\left(\eta,\zeta\right)^{1/2}} \Psi_n\left(x_0 \et\right) \Psi_n\left(\sigma\right),
		\end{equation*}
		which readily gives
		\begin{equation*}
			\Hat{E}_3\left(t,x,\eta,\zeta;x_0\right)=-H\left(t\right)\et\sum_{n=0}^\infty\frac{\sin\left(\rho_n\left(\eta,\zeta\right)^{1/2} t\right)}{\rho_n\left(\eta,\zeta\right)^{1/2}} \Psi_n\left(x_0 \et\right) \Psi_n\left(\et x\right)
		\end{equation*}
		so that
		\begin{align*}
			E_3&\left(t,x,y,z;x_0\right)=-\frac{H\left(t\right)}{2\pi}\\\times&\int_{\bb{R}^2}e^{iy\eta+iz\zeta}\et\sum_{n=0}^\infty\frac{\sin\left(\rho_n\left(\eta,\zeta\right)^{1/2} t\right)}{\rho_n\left(\eta,\zeta\right)^{1/2}}\Psi_n\left(x_0 \et\right) \Psi_n\left(x \et\right) d\eta\, d\zeta.
		\end{align*}
		Now, we are interested in 
		\[\mathtt{F}_3\left(t-s,x,y,z;\xi_0,\heta,\hzeta\right):=\cl{F}E\left(t-s,x,y-\cdot,z-\cdot,\cdot\right)\left(\xi_0,\heta,\hzeta\right);\]
		by translation and dilation properties of the Fourier transform, we see
		\begin{align*}
			\mathtt{F}_3\left(t-s,x,y,z;\xi_0,\heta,\hzeta\right)=e^{iy\heta+iz\hzeta}\int_{-\infty}^\infty e^{-i\xi_0 x_0}\Hat{E_3}\left(t-s,x,-\heta,-\hzeta;x_0\right). 
		\end{align*}
		Thus, we finally get
		\begin{align}
			\label{F}
			\mathtt{F}_3&\left(t-s,x,y,z;\xi_0,\heta,\hzeta\right)=-\frac{H\left(t-s\right)}{2\pi}e^{iy\heta+iz\hzeta}\nonumber\\\times&\sum_{n=0}^\infty\left(-i\right)^n\frac{\sin\left(\rho_n\left(-\heta,\hzeta\right)^{1/2} \left(t-s\right)\right)}{\rho_n\left(-\heta,\hzeta\right)^{1/2}} \Psi_n\left(\frac{\xi_0}{\het}\right) \Psi_n\left(\het x\right).
		\end{align}
		Following the above procedure along, one can also solve
		\begin{equation*}
			P_2E_2=\delta(t)\delta(x-x_0)\delta(y)
		\end{equation*}
		and obtain
		\begin{align*}
			E_2&\left(t,x,y;x_0\right)=-\frac{H\left(t\right)}{2\pi}\\\times&\int_{\bb{R}}e^{iy\eta}\et\sum_{n=0}^\infty\frac{\sin\left(\hat{\rho}_n\left(\eta\right)^{1/2} t\right)}{\hat{\rho}_n\left(\eta\right)^{1/2}}\Psi_n\left(x_0 \et\right) \Psi_n\left(x \et\right) d\eta,
		\end{align*}
		where $\hat{\rho}_n:=\rho_n(\eta,0)$. And thus
		\begin{align}
			\label{F2}
			\mathtt{F}_2&\left(t-s,x,y;\xi_0,\heta\right)=-\frac{H\left(t-s\right)}{2\pi}e^{iy\heta}\nonumber\\\times&\sum_{n=0}^\infty\left(-i\right)^n\frac{\sin\left(\hat{\rho}_n\left(-\heta)\right)^{1/2} \left(t-s\right)\right)}{\hat{\rho}_n\left(-\heta\right)^{1/2}} \Psi_n\left(\frac{\xi_0}{\het}\right) \Psi_n\left(\het x\right).
		\end{align}      
\subsection{Two dimensions: Integral bound}
\label{two dim}
Again  we rewrite here the expression for the two-dimensional candidate random field solution
	\begin{equation*}
	u_2(t,x,y)=\int_0^t \int_{\bb{R}^d}E_2\left(t-s,x,y-y_0;x_0\right)\Dot{F}_1(s,y_0,x_0)ds\,dx_0\,dy_0.
\end{equation*}
In the present section, we proceed with the study of the fundamental solutions $E_2$ to argue that the formal expression above is well defined as a real-valued process. We start proving
\begin{lemma}
	\label{lemma 1 2 dim}
	Let $\mathtt{F}_2\left(t-s,x,y;\xi_0,\heta\right)$ be defined as in (\ref{F2}) and the assumptions of Theorem \ref{theorem 1} (ii) hold. Then 
	\begin{equation*}
		I_2:=\int_0^t \int_{\bb{R}^2}\left|\mathtt{F}_2\left(t-s,x,y;\xi_0,\heta\right)\right|^2 \nu_2\left(d\xi_0,d\heta\right)\,ds <+\infty.
	\end{equation*}
\end{lemma}
\begin{proof}
	Let us call
	\begin{equation*}
		J\left(t-s,x,y;\heta\right):=\int_{\bb{R}}\left|\mathtt{F}_2\left(t-s,x,y;\xi_0,\heta\right)\right|^2 \,d\xi_0
	\end{equation*}
	and change variable in the following way: $\lambda=|\heta|^{-1/2}\xi_0$. Then
	\begin{equation*}
		J\left(t-s,x,y;\heta\right):=\int_{\bb{R}}\left|\mathtt{F}_2\left(t-s,x,y;|\heta|^{1/2}\lambda,\heta\right)\right|^2 |\heta|^{1/2}\,d\lambda.
	\end{equation*}
	Now, if we denote
	\begin{align}
		\label{g_n}
		g_n\left(t-s,x,y;\heta\right):=&\frac{H\left(t-s\right)}{2\pi}e^{iy\heta}\nonumber\\\times&\left(-i\right)^n\frac{\sin\left(\hat{\rho}_n\left(-\heta\right)^{1/2} \left(t-s\right)\right)}{\hat{\rho}_n\left(-\heta\right)^{1/2}} \Psi_n\left(\het x\right),
	\end{align}
	we get
	\begin{equation*}
		\mathtt{F}_2\left(t-s,x,y;|\heta|^{1/2}\lambda,\heta\right)=\sum_{n=0}^\infty g_n\left(t-s,x,y;\heta\right) \Psi_n (\lambda)
	\end{equation*}
	and hence, by Parseval's identity, 
	\begin{align}
		\label{parseval}
		\int_{\bb{R}}\left|\mathtt{F}_2\left(t-s,x,y;|\heta|^{1/2}\lambda,\heta\right)\right|^2 \,d\lambda=\sum_{n=0}^\infty g_n\left(t-s,x,y;\heta\right)^2.
	\end{align}
	It follows that
	\begin{align*}
		I_2=&\int_0^t \int_{\bb{R}}|\heta|^{1/2}\sum_{n=0}^\infty g_n\left(t-s,x,y;\heta\right)^2\hat{\nu}_2\left(d\heta\right)\,ds\\=&\frac{1}{4\pi^2}\int_0^t \int_{\bb{R}}\sum_{n=0}^\infty\frac{\sin^2\left(\hat{\rho}_n^{1/2}(\heta)\left(t-s\right)\right)}{\hat{\rho}_n(\heta)}|\heta|^{1/2}\Psi_n^2\left(|\heta|^{1/2}x\right)\hat{\nu}_2\left(d\heta\right)\,ds.
	\end{align*}
	Recalling (\ref{bound on psi2}) and denoting 
	\begin{equation*}
		T_0:=\int_0^t \int_{\bb{R}} g_0\left(t-s,x,y;\heta\right)^2|\heta|^{1/2}\hat{\nu}_2\left(d\heta\right)\,ds,
	\end{equation*}
	we get
	\begin{align}
		\label{I bound}
		I_2 \leq &\frac{D}{4\pi^2}\int_0^t \int_{\bb{R}^2}\sum_{n=1}^\infty\frac{|\heta|^{1/2}}{n^{1/6}\left(\left(2n+1\right)\mu\left|\heta\right|+b\heta\right)}\hat{\nu}_2\left(d\heta\right)\,ds+T_0\nonumber\\\leq &\frac{D\,t}{4\pi^2n^{1/6}}\sum_{n=1}^\infty\int_{\bb{R}^2}\frac{|\heta|^{1/2}}{\left(\left(2n+1\right)\mu\left|\heta\right|+b\heta\right)}\hat{\nu}_2\left(d\heta\right)+T_0.
	\end{align}
	Now, we have
	\begin{align*}
		 \int_{\bb{R}}\frac{|\heta|^{1/2}}{\left(\left(2n+1\right)\mu\left|\heta\right|+b\heta\right)}\hat{\nu}_2\left(d\heta\right)\leq\frac{1}{n\mu}\int_{\bb{R}} \left|\heta\right|^{-1/2}\hat{\nu}(d\heta)=\frac{C}{\mu n}
	\end{align*}
	where the last inequality comes from the fact that the condition $|b|< \mu $ gives $\mu  + \sgn(\heta)b >0$. 
	Finally, convergence of the integral $T_0$ under this measure is easily checked. Indeed, we can find an $\epsilon>0$ s.t. $\mu+\sgn(\heta)b>\epsilon$ for all $\heta \in \bb{R}$ and observe that $\Psi_0^2\left(|\heta|^{1/2}x\right)=\pi^{-1/4}e^{-\left|\heta\right|x^2/2}\leq \pi^{-1/4}$. Therefore
	\begin{align*}
		T_0 \leq  \frac{t}{4\pi^{9/4}\epsilon} \int_{\bb{R}} \left|\heta\right|^{-1/2} \hat{\nu}_2\left(d\heta\right)<+\infty
	\end{align*}
\end{proof}
Now, we want to find a sequence $E_{2,m}\left(t,x,y;x_0\right) \in \cl{E}_0$ so that $E_{2,m} \to E_2$ in $\Bar{\cl{P}}$. We set
\begin{align*}
	E_{2,m}\left(t,x,y;x_0\right):=&-\frac{H\left(t\right)}{2\pi}\\\times&\int_{\bb{R}}e^{iy\eta}\et\sum_{n=0}^m\frac{\sin\left(\hat{\rho}_n\left(\eta\right)^{1/2} t\right)}{\hat{\rho}_n\left(\eta\right)^{1/2}}\Psi_n\left(x_0 \et\right) \Psi_n\left(x \et\right) d\eta.
\end{align*}
Again for each $m \in \bb{N}$, $E_m \in \cl{S}\left(\bb{R}^2\right)$, and so calling
\begin{align*}
	\mathtt{F}_2^m\left(t-s,x,y;\xi_0,\heta\right):=&\frac{H\left(t-s\right)}{2\pi}e^{iy\heta}\\\times&\sum_{n=m}^{+\infty}\left(-i\right)^n\frac{\sin\left(\hat{\rho}_n\left(-\heta\right)^{1/2} \left(t-s\right)\right)}{\hat{\rho}_n\left(-\heta\right)^{1/2}} \Psi_n\left(\frac{\xi_0}{\het}\right) \Psi_n\left(\het x\right),
\end{align*}
we have
\begin{align*}
	\left|\left|E_{2,m}-E_2\right|\right|_0=&\int_0^t ds\int_{\bb{R}^2}\left|\mathtt{F}_2^m\left(t-s,x,y;\xi_0,\heta\right)\right|^2\nu_2\left(d\xi_0,d\heta\right)\\\leq &\frac{D\,t}{4\pi^2}\sum_{n=m}^{+\infty}\int_{\bb{R}^2}\frac{|\heta|^{1/2}}{n^{1/6}\left(\left(2n+1\right)\mu\left|\heta\right|+b\heta\right)}\left(1+\heta^2\right)^{-\alpha}d\heta  \\ \leq &\Tilde{D}\,t\sum_{n=m}^{+\infty} n^{-7/6} \xrightarrow{m \to +\infty}0.
\end{align*}

\subsection{Three dimensions: integral bound}
\label{3dib}
For the reader's convenience, we recall the formal expression for the candidate three-dimensional random field solution 
\begin{align*}
	u_3(t,x,y,z)=\int_0^t \int_{\bb{R}^d}&E_3\left(t-s,x,y-y_0,z-z_0;x_0\right)\\ \times&\Dot{F}_3(s,x_0,y_0,z_0)ds\,dx_0\,dy_0\,dz_0.
\end{align*}
In the present section, we would like to understand whether the above epxression is a real-valued process and to do that we need to check if $E_3 $ belongs to $ \cl{P}_0$. As always, we start with showing its membership to $\Bar{\cl{P}}$.
\begin{lemma}
	\label{lemma 1}
	Let $\mathtt{F}_3\left(t-s,x,y,z;\xi_0,\heta,\hzeta\right)$ be defined as in \ref{F} and the assumptions of Theorem \ref{theorem 1} (iii) hold. Then 
	\begin{equation*}
		I_3:=\int_0^t \int_{\bb{R}^3}\left|\mathtt{F}_3\left(t-s,x,y,z;\xi_0,\heta,\hzeta\right)\right|^2 \nu_3\left(d\xi_0,d\heta,d\hzeta\right)\,ds <\infty.
	\end{equation*}
\end{lemma}
\begin{proof}
	The proof is very similar to the one of Proposition \ref{lemma 1 2 dim}. Indeed calling
	\begin{align*}
		\hat{g}_n\left(t-s,x,y;\heta,\hzeta\right):=&\frac{H\left(t-s\right)}{2\pi}e^{iy\heta+iz\hzeta}\nonumber\\\times&\left(-i\right)^n\frac{\sin\left(\rho_n\left(-\heta,\hzeta\right)^{1/2} \left(t-s\right)\right)}{\rho_n\left(-\heta,\hzeta\right)^{1/2}} \Psi_n\left(\het x\right),
	\end{align*}
	and
	\begin{gather*}
		\hat{T}_0:=\int_0^t \int_{\bb{R}^2} \hat{g}_0\left(t-s,x,y,z;\heta,\hzeta\right)^2|\heta|^{1/2}\hat{\nu}_3\left(d\heta,d\hzeta\right)\,ds,\\	
		\hat{q}_n:= \int_{\bb{R}^2}\frac{|\heta|^{1/2}}{\left(\left(2n+1\right)\mu\left|\heta\right|+a\hzeta^2+b\heta\right)}\hat{\nu}_3\left(d\heta,d\hzeta\right),
	\end{gather*}
	with the same strategy, we obtain
	\begin{align}
		\label{Ihat bound}
		I_3 \leq \frac{D\,t}{4\pi^2}\sum_{n=1}^\infty n^{-1/6} q_n+\hat{T}_0.
	\end{align}
	Now, since $|b|< \mu $, we have $\mu  + \sgn(\mu)b >0$ and consequently
	\begin{equation*}
		q_n \leq \int_{\bb{R}^2} \frac{|\heta|^{1/2}}{n\mu\left|\heta\right|+a\hzeta^2}w(\heta^2+\hzeta^2)\,d\heta d\hzeta=:J(n).
	\end{equation*}
	We can then switch to polar coordinates and utilize Young's inequality 
	\begin{equation*}
		\frac{x^p}{p}+\frac{y^q}{q} \geq xy, \quad p+q=1,
	\end{equation*}
	to get
	\begin{align*}
		J(n)=&4\int_0^{+\infty} d\rho\; \rho\, w(\rho^2)\int_0^{\pi/2} \frac{\left(\rho \cos(\theta)\right)^{1/2}}{n \mu \rho \cos(\theta) + a \rho^2 \sin^2(\theta)}d\theta\\= &4\int_0^{+\infty}d\rho\;\rho^{1/2}w(\rho^2) \int_0^{\pi/2} \frac{ \left(\cos(\theta)\right)^{1/2}}{n \mu \cos(\theta) + a \rho \sin^2(\theta)}d\theta\\\leq&\frac{4p^pq^q}{\mu^pa^q}n^{-p}\int_0^{+\infty}d\rho\;\rho^{1/2-q}w(\rho^2) \int_0^{\pi/2} \cos(\theta)^{1/2-p}\sin{\theta}^{-2q}d\theta.
	\end{align*}
	We may very well select $0<\delta<1/6$ and $5/6+\delta>p>5/6$ from which it follows that $1/6-\delta<q<1/6$. This choice is sufficient for the convergence of the integral in $d\theta$. Indeed, we get
	\begin{equation*}
		\int_0^{\pi/2} \cos(\theta)^{1/2-p}\sin{\theta}^{-2q}d\theta=\frac{\G\left(3/4-p/2\right)\G\left(1/2-q\right)}{2\G\left(5/4-p/2-q\right)}<+\infty.
	\end{equation*}
	Furthermore, condition (\ref{nu3 condition}) guarantees we can select $\delta$ so that
	\begin{equation*}
		\int_0^{+\infty}\rho^{1/2-q}w(\rho^2)d\rho < \int_0^{+\infty}\rho^{1/3+\delta}w(\rho^2)d\rho <+\infty.
	\end{equation*}
	Therefore, we get $J(n)\leq\hat{C}n^{-p}, \ p>5/6$ which is sufficent for the convergence of the series in (\ref{Ihat bound}).\\
	Finally, taking $\epsilon>0$ s.t. $\epsilon<\mu+b\sgn(\heta))$ we write
	\begin{align*}
		T_0 \leq &\frac{t}{4\pi^{9/4}}\int_{\bb{R}^2} \frac{|\heta|^{1/2}}{a\hzeta^2+\epsilon\left|\heta\right|} w(\heta^2+\hzeta^2)\,d\heta d\hzeta,
	\end{align*}
	which is finite for the same arguments utilized above.
\end{proof}

Having proved that $E_3 \in \Bar{\cl{P}}$, in order to show that $E_3 \in \cl{P}_0$, we need to find a sequence of functions $E_{3,m} \in \cl{E}_0$ such that $E_{3,m} \to E_3$ in $\Bar{\cl{P}}$. We set
\begin{align*}
	E_{3,m}\left(t,x,y,z;x_0\right):=&-\frac{H\left(t\right)}{2\pi}\\\times&\int_{\bb{R}^2}e^{iy\eta+iz\zeta}\et\sum_{n=0}^m\frac{\sin\left(\rho_n\left(\eta,\zeta\right)^{1/2} t\right)}{\rho_n\left(\eta,\zeta\right)^{1/2}}\Psi_n\left(x_0 \et\right) \Psi_n\left(x \et\right) d\eta \, d\zeta.
\end{align*}
It is readily seen that $E_{3,m} \in \cl{S}\left(\bb{R}^3\right)$ for all $m \in \bb{N}$, 
and
\begin{align*}
	\left|\left|E_{3,m}-E\right|\right|_0 \leq &\frac{D\,t}{4\pi^2}\sum_{n=m}^{+\infty}n^{-1/6}\hat{q}_n \xrightarrow{m \to +\infty}0
\end{align*}

\subsection{Common computations: measurability}
\label{fmeas}
To prove our claim, it is left to show that the maps $(t,x,y)\mapsto u_2(t,x,y)$ and $(t,x,y,z)\mapsto u_3(t,x,y,z)$ are measurable. Our strategy in this regard is to show $L_2(\Omega)$-continuity of the solution processes.
We only tackle the three-dimensional case, as the two-dimensional case is pretty much analogous.

We first consider the variable $y$. We call
\begin{align*}
	\Bar{F}_3(t-s,x,z;\xi_0,\heta,\hzeta):=&\frac{H\left(t-s\right)}{2\pi}e^{iz\hzeta}\\\times&\sum_{n=0}^{+\infty}\left(-i\right)^n\frac{\sin\left(\rho_n\left(-\heta,\hzeta\right)^{1/2} \left(t-s\right)\right)}{\rho_n\left(-\heta,\hzeta\right)^{1/2}} \Psi_n\left(\frac{\xi_0}{\het}\right) \Psi_n\left(\het x\right),
\end{align*}
and get
\begin{align*}
	&\bb{E}\left[\left|u_3(t,x,y+h,z)-u_3(t,x,y,z)\right|^2\right]\\=&\int_0^tds\int_{\bb{R}^3} \left|\mathtt{F}_3\left(t-s,x,y+h,z;\xi_0,\heta,\hzeta\right)- \mathtt{F}_3\left(t-s,x,y,z;\xi_0,\heta,\hzeta\right)\right|^2\nu_3(d\xi_0,d\heta,d\hzeta)\\=&\int_0^tds\int_{\bb{R}^3} \left|e^{i\heta(y+h)}-e^{i\heta y}\right|^2\left|\Bar{F}_3(t-s,x,z;\xi_0,\heta,\hzeta)\right|^2\nu_3(d\xi_0,d\heta,d\hzeta).
\end{align*}
Now, $\Bar{F}(t-s,x;\xi_0,\heta)$ is integrable and, hence, by dominated convergence, we have
\begin{equation*}
	\lim_{h\to 0}\bb{E}\left[\left|u_3(t,x,y+h,z)-u_3(t,x,y,z)\right|^2\right]=0.
\end{equation*}
The $z$ variable is treated the same.\\
For the increment in the $x$ variable, we denote
\begin{equation*}
	S_3(t-s,\heta,\hzeta,x,h):=\sum_{n=0}^{+\infty}\frac{\sin^2\left(\rho_n\left(-\heta,\hzeta\right)^{1/2}(t-s)\right)}{\rho_n\left(-\heta,\hzeta\right)}\left(\Psi_n(\et(x+h))-\Psi_n(\et x)\right)^2,
\end{equation*} 
and write
\begin{align*}
	&\bb{E}\left[\left|u_3(t,x+h,y,z)-u_3(t,x,y,z)\right|^2\right]\\=&\int_0^tds\int_{\bb{R}^3} \left|\mathtt{F}_3\left(t-s,x+h,y,z;\xi_0,\heta,\hzeta\right)- \mathtt{F}_3\left(t-s,x,y,z;\xi_0,\heta,\hzeta\right)\right|^2\nu_3(d\xi_0,d\heta,d\hzeta)\\=&\int_0^tds\int_{\bb{R}^2} \frac{|\heta|^{1/2}}{4\pi^2}S_3(t-s,\heta,\hzeta,x,h)\left(1+|\heta|^2+|\hzeta|^2\right)^{-\alpha}d\heta d\hzeta.
\end{align*}
Moreover thanks to (\ref{bound on psi2}), we have
\begin{align}
	\label{sigma bound}
	S_3(t-s,\heta,\hzeta,x,h) \leq \sum_{n=0}^{+\infty} \frac{4D}{n^{1/6}\rho_n\left(-\heta,\hzeta\right)}< +\infty
\end{align}
granting uniform convergence and consequently the continuity of $S_3(t-s,\heta,\hzeta,x,h)$ w.r.t. his last variable. Furthermore, in the proof of Lemma \ref{lemma 1} we have seen  that
\begin{align*}
	\int_0^tds\int_{\bb{R}^3} \et \sum_{n=0}^{+\infty} \frac{1}{n^{1/6}\rho_n\left(-\heta,\hzeta\right)}\nu_3(d\xi_0,d\heta,d\hzeta)< + \infty.
\end{align*}
Therefore, by dominated convergence we get
\begin{equation}
	\lim_{h\to 0}\bb{E}\left[\left|u_3(t,x+h,y,z)-u_3(t,x,y,z)\right|^2\right]=0.
\end{equation}
Lastly, the increment in time yields
\begin{align*}
	&\bb{E}\left[\left|u_3(t+h,x,y,z)-u_3(t,x,y,z)\right|^2\right]\\ \leq&2\int_0^tds\int_{\bb{R}^3} \left|\mathtt{F}_3\left(t+h-s,x,y,z;\xi_0,\heta,\hzeta\right)- \mathtt{F}_3\left(t-s,x,y,z;\xi_0,\heta,\hzeta\right)\right|^2\nu_3(d\xi_0,d\heta,d\hzeta) \\\mathrel{\phantom{=}}&+2\int_t^{t+h}ds\int_{\bb{R}^3} \left|\mathtt{F}_3\left(t+h-s,x,y,z;\xi_0,\heta,\hzeta\right)\right|^2\nu_3(d\xi_0,d\heta,d\hzeta)\\=&2J_1(h)+2J_2(h).
\end{align*}
Lemma \ref{lemma 1} immediately gives
\begin{equation*}
	\lim_{h \to 0}J_2(h)=0. 
\end{equation*}
Denoting
\begin{equation*}
	S_3'(t-s,\heta,\hzeta, x,h):=\sum_{n=0}^{+\infty}\frac{\left(\sin\left(\rho_n^{1/2}(\heta,\hzeta)(t+h-s)\right)-\sin\left(\rho_n^{1/2}(\heta,\hzeta)(t-s)\right)\right)^2}{\rho_n(\heta,\hzeta)}\Psi_n^2(\et x),
\end{equation*}
by Parseval's identity, we obtain
\begin{align*}
	J_1(h)=\int_0^tds\int_{\bb{R}^2} \frac{|\heta|^{1/2}}{4\pi^2}S_3'(t-s,\heta,\hzeta, x,h)\left(1+|\heta|^2+|\hzeta|^2\right)^{-\alpha}d\heta d\hzeta.
\end{align*}
For $S_3'(t-s,\heta,\hzeta, x,h)$, the same bound as in (\ref{sigma bound}) holds. So, by dominated convergence, again we have
\begin{equation*}
	\lim_{h\to 0}J_1(h)=0.
\end{equation*}



\nocite{*}
\printbibliography
\end{document}